\documentclass[11pt,a4paper]{article}
\usepackage{amsmath}
\usepackage{amssymb}
\usepackage{latexsym}
\newcommand{\gM}{{\mathfrak M}}

\newcommand{\dual}{\makebox[0mm]{}^{{\scriptstyle\vee}}}
\newtheorem{theorem}{Theorem}[section]

\newtheorem{lemma}[theorem]{Lemma}
\newtheorem{proposition}[theorem]{Proposition}

\newtheorem{corollary}[theorem]{Corollary}

\newtheorem{exmp}[theorem]{Example}
\newtheorem{exmps}[theorem]{Examples}
\newtheorem{rem}[theorem]{Remark}

\newcommand{\prf}{{\em Proof}. }
\newcommand{\qed}{\hspace*{\fill}$\Box$}

\newcommand{\beeq}[1]{\begin{eqnarray}\label{#1}}
\newcommand{\eneq}{\end{eqnarray}}

\newcommand{\kc}{{\cal C}}

\newcommand{\kp}{{\cal P}}

\newcommand{\IC}{{\mathbb C}}

\newcommand{\IP}{{\mathbb P}}
\newcommand{\IQ}{{\mathbb Q}}
\newcommand{\IR}{{\mathbb R}}
\newcommand{\IZ}{{\mathbb Z}}

\newcommand{\Pic}{{\rm Pic}}

\newcommand{\Aut}{{\rm Aut}}
\newcommand{\vol}{{\rm vol}}
\newcommand{\SSp}{{\rm Sp}}
\newcommand{\td}{{\rm td}}

\newcommand{\verylongarrow}[1]{\hbox to #1{\rightarrowfill}}

\begin{document}
{\centerline{\Large\bf  Finiteness results for compact
hyperk\"ahler manifolds}}

\bigskip

{\centerline{\large Daniel Huybrechts}

\bigskip

\bigskip

\bigskip

\bigskip

In \cite{S} Sullivan proved the following finiteness result for
simply connected K\"ahler manifolds.

\bigskip

{\bf Theorem} ---
{\it The diffeomorphism type of a simply connected K\"ahler manifold is 
determined up to finite ambiguity by its integral cohomology ring and
its Pontrjagin classes.}

\bigskip

The aim of this note is to point out that much less of the topology needs
to be fixed in order to determine the diffeomorphism type of a compact
hyperk\"ahler manifold (up to finite possibilities).
Recall that a hyperk\"ahler manifold is a $4n$-dimensional Riemannian
manifold $(M,g)$ with holonomy $\SSp(n)$. If such a manifold
is compact it is automatically
simply connected \cite[Prop. 4]{B}.
The first Pontrjagin class $p_1(M)\in H^4(M,\IZ)$ defines 
a homogeneous polynomial of degree
$2n-2$ on $H^2(M,\IZ)$ by $\alpha\mapsto\int_M\alpha^{2n-2}p_1(M)$.
\bigskip

{\bf Theorem \ref{mainthm}} --- {\it If the second integral cohomology $H^2$
and the homogeneous polynomial of degree $2n-2$ on $H^2$ defined
by the first Pontrjagin class are given, then there exist at most
finitely many diffeomorphism types of compact
hyperk\"ahler manifolds of dimension $4n$ realizing this structure.}

\bigskip

 We will also see
that instead of fixing the action of $p_1$ in the above theorem one
could as well fix $H^2$ with the top intersection product on $H^2$.
This would amount to fixing the action of $p_0(M)$.

The proof of this theorem uses the finiteness result of Koll\'ar and
Matsusaka \cite{KM}, a formula  by Hitchin and Sawon \cite{HS}, the
projectivity criterion for hyperk\"ahler manifolds, and the surjectivity
of the period map proved in \cite{H}.

On any compact hyperk\"ahler manifold $(M,g)$ there exists a $\IP^1$ of complex
structures on $M$ compatible with the metric $g$. With any of these complex
structures the manifold $M$ becomes an irreducible holomorphic
symplectic manifold. Conversely, any K\"ahler class on an irreducible
holomorphic symplectic manifold is uniquely induced by a hyperk\"ahler metric
on the underlying real manifold $M$. Once the diffeomorphism
type of $M$ is fixed one wants to know how many deformation
types of hyperk\"ahler
metrics $g$ or, equivalently, of irreducible holomorphic symplectic complex
structures do exist on $M$. Again, there is only a finite number
of possibilities.

\bigskip

{\bf Theorem \ref{mainthm2}} --- {\it Let $M$ be a fixed compact manifold. Then
there exist at most finitely many different deformation types of irreducible
holomorphic symplectic complex  structures on $M$.}


\section{Finiteness of polarized manifolds}

Let $X$ be a projective manifold over $\IC$. If $L\in\Pic(X)$ is an ample
line bundle on $X$, then the Hilbert polynomial of the polarized
manifold $(X,L)$ is $P(X,L,z)\in\IQ[z]$, such that
$P(X,L,m)=\chi(X,L^{\otimes m})$ for all $m\in\IZ$.
For some high power $m$ the line bundle
$L^{\otimes m}$ is very ample and $X$ is naturally embedded into a projective
space $\IP^N$ of dimension $N=P(X,L,m)$. The Quot-scheme parametrizing
all subvarieties of $\IP^N$ with given Hilbert polynomial is projective
and therefore has only finitely many irreducible components. In particular,
the number of deformation types of smooth projective varieties parametrized
by this Quot-scheme is finite. In \cite{M} Matsusaka proved
that there exists a universal power $m$ as above that depends only
on the numerical polynomial $P(X,L,z)\in\IQ[z]$ and not on $X$ itself.

On the other hand, for a given $d,a_0,a_1\in\IQ$ there exists only a finite
number of numerical polynomials $P_i(z)\in\IQ[z]$ of degree
$d$ with $P_i(z)=a_0z^d+a_1z^{d-1}+\ldots$ that occur as Hilbert polynomials
of polarized manifolds $(X,L)$ as above \cite{KM}.

Note that
$P(X,L,z)=\frac{L^d}{d!}z^d-\frac{(K_X.L^{d-1})}{2(d-1)!}z^{d-1}+\ldots$. Thus
one obtains the following result due to Koll\'ar and Matsusaka:

\begin{theorem}{\bf \cite[Thm.\ 3]{KM}} ---
The number of deformation types of projective manifolds $X$ of dimension $d$
that admit an ample line bundle $L$ with fixed intersections
$L^d, (K_X.L^{d-1})\in\IZ$ is finite.
\end{theorem}

For manifolds with trivial canonical bundle $K_X$ this becomes

\begin{corollary}\label{CorKM}---
The number of deformation types of projective manifolds $X$ of dimension
$d$ with trivial canonical bundle that admit an ample line bundle $L$ with
bounded $L^d\in\IZ$ is finite.
\end{corollary}

The result in particular applies to irreducible holomorphic symplectic
manifolds. In order to get rid of the supplementary polarization
of the manifold one has to show that any irreducible holomorphic symplectic
manifold can be deformed to one that admits a polarization $L$
with bounded $L^d$. This will be done in the next section.


\section{Complex structures on a fixed manifold}\label{complex}

Here we will prove the following

\begin{theorem}\label{mainthm2}---
Let $M$ be a fixed compact manifold. Then there exist at most finitely many
different deformation types of irreducible holomorphic symplectic
complex  structures on $M$.
\end{theorem}

Recall that an irreducible holomorphic symplectic manifold is a compact
complex K\"ahler manifold $X$, which is simply connected and admits
a nowhere degenerate two-form $\sigma$ such that
the space of global holomorphic two-forms $H^0(X,\Omega^2_X)$ is generated by
$\sigma_X$. It will be clear from the proof
that `deformation type' could also be taken in the more restrictive
sense that just allows families of irreducible holomorphic symplectic
manifolds.

If $I$ is a complex structure on $M$ such that $X=(M,I)$ is irreducible
holomorphic symplectic, then the Beauville-Bogomolov form $q_X$ is a primitive
integral quadratic form on $H^2(M,\IZ)$ with the property that for some
positive constant $c$ $$q_X(\alpha)^n=c\int_X\alpha^{2n}.\eqno{(1)}$$ Clearly,
the integral form $q_X$ does not change under deformation of $I$. But in fact
$q_X$ is completely independent of the complex structure. The argument to prove
this combines results of Fujiki and Nieper and will be given in Sect.\
\ref{unnorm}. We will denote this quadratic form on $H^2(M,\IZ)$ by $q$. Thus,
whenever $X=(M,I)$ is irreducible holomorphic symplectic, then
$(H^2(M,\IZ),q)\cong (H^2(X,\IZ),q_X)$.

Let $\Gamma$ be the lattice $(H^2(M,\IZ),q)$. The moduli space of
marked irreducible holomorphic symplectic manifolds $\gM_\Gamma$
is defined
as $\gM_\Gamma:=\{(X,\varphi)|\varphi:(H^2(X,\IZ),q_X)\cong\Gamma\}/\sim$.
Here, $\sim$ is the equivalence relation induced by $\pm f^*$, where
$f:X\to X'$ is an arbitrary biholomorphic map.
One proves that $\gM_\Gamma$ has a natural smooth complex structure, which
however is not Hausdorff.
The period map is the holomorphic map
$\kp:\gM_\Gamma\to Q_\Gamma\subset\IP(\Gamma\otimes\IC)$,
$(X,\varphi)\mapsto \varphi(H^{2,0}(X))$, where $Q_\Gamma$ is the period
domain $Q_\Gamma:=\{x|q(x)=0,~q(x+\bar x)>0\}$. We will
use the following result \cite[Thm.\ 8.1]{H}

\begin{theorem}\label{surj}--- If $\gM_\Gamma^o$ is a non-empty
connected component of
$\gM_\Gamma$, then $\kp:\gM_\Gamma^o\to Q_\Gamma$ is surjective.
\end{theorem}

Let us fix a primitive element $0\ne\alpha\in\Gamma$ with $q(\alpha)>0$ and
consider the intersection $S_\alpha$ of $Q_\Gamma$ with the hyperplane
$H_\alpha$ defined by $q(x,\alpha)=0$. The assumption that $\alpha$ is
primitive is not important, but makes the formulation of some statements
more elegant.

\begin{lemma}\label{notempty}---
For any $0\ne\alpha\in\Gamma$ the intersection $S_\alpha$ is not empty.
\end{lemma}

\prf Let $V$ be a real vector space with a non-degenerate quadratic form
$(~,~)$ of index $(p,q)$.
If $p\geq2$ then one finds orthogonal vectors
$y_1,y_2\in V$, such that $(y_1,y_1)=(y_2,y_2)>0$. Then
$y=y_1+i y_2\in V\otimes\IC$ satisfies $(y,y)=0$ and $(y+\bar y,y+\bar y)>0$.
Apply this to the orthogonal complement $V$ of $\alpha\in\Gamma$
in $\Gamma\otimes_\IZ\IR$. Since $q$ has index $(3,b_2(M)-3)$, the induced
quadratic form on $V$ has two positive eigenvalues.\qed

\begin{lemma}--- If $y\in S_\alpha$ is generic, then
$\{\beta\in\Gamma|q(y,\beta)=0\}=\IZ\alpha$.
\end{lemma}

\prf Clearly, $\IZ\alpha\subset\{\beta\in\Gamma|q(y,\beta)=0\}$.
If for generic $y$ the right hand side were strictly bigger, then
$S_\alpha$ would equal $S_\beta$ for some $\beta\in\Gamma$
linearly independent of $\alpha$.
Thus, the two hyperplanes defined by $\alpha$ and $\beta$, resp., would
contain the intersection of the quadric defined by $q$ and $H_\alpha$.
Hence, $q$ restricted to $H_\alpha$ would be degenerate, which contradicts
$q(\alpha)>0$.\qed

\begin{proposition}---
If $(Y,\psi)\in\gM_\Gamma^o$ such that $y:=\kp(Y,\psi)\in S_\alpha$
is generic, then $\varphi:\Pic(Y)\cong\IZ\alpha$ and $\Pic(Y)$
is generated by an ample line bundle $L$. In particular, $Y$ is projective.
\end{proposition}

\prf One knows that $\Pic(Y)=\{\beta\in H^2(Y,\IZ)|q_Y(\beta,\sigma_Y)=0\}$.
By the previous lemma, the first assertion follows from this. Thus,
$\Pic(Y)$ is generated by a line bundle $L$ with $q_Y(L)>0$. Hence,
the intersection of the positive cone $\kc_Y$ with $H^2(Y,\IZ)$
is non-empty and by the projectivity criterion \cite{H} this shows that $Y$
is projective. Therefore, $L$ or $L\dual$ is ample.\qed

\bigskip

{\bf Remark} --- The proof of the projectivity criterion in \cite{H} was
flawed. A complete proof was given recently in the Erratum \cite{HErr} using a
new result by Demailly and Paun \cite{DP}.

\bigskip

Thus, for any fixed primitive $0\ne\alpha\in\Gamma$ with $q(\alpha)>0$ one
can deform any $(X,\varphi)\in\gM_\Gamma$ to a marked manifold
$(Y,\psi)\in\gM_\Gamma$ (which then is of course contained
in the same connected component of $\gM_\Gamma$) such that
$(Y,L:=\psi^{-1}(\alpha))$ is a polarized projective manifold
with fixed $q_Y(L)=q(\alpha)$ (use Thm.\ \ref{surj} and Lemma
\ref{notempty}). In order to apply the boundedness result
of Koll\'ar and Matsusaka one has to fix $L^{2n}$ instead.
Since the underlying real manifold for the different components of
$\gM_\Gamma$ and thus the constant $c$ in (1)
could be different, fixing $q_Y(L)$ is a priori not enough in order
to fix $L^{2n}$.

\bigskip

{\it Proof of Thm.\ \ref{mainthm2}.} Let $\alpha\in H^2(M,\IZ)$ be fixed such that
$\alpha$ is primitive and $q(\alpha)>0$.
Let $I$ be a complex structure on $M$, such that $X=(M,I)$ is an irreducible
holomorphic symplectic manifold.
Then fixing a marking $\varphi:(H^2(X,\IZ),q_X)\cong
\Gamma\cong (H^2(M,\IZ),q)$ yields a point $(X,\varphi)$ in the moduli
space $\gM_\Gamma$. Hence, $X$ is deformation equivalent to a projective
manifold $Y$ that admits an ample line bundle $L$ with $L^{2n}=c^{-1}q_Y(L)=
c^{-1}q(\alpha)$. By Koll\'ar and Matsusaka (cf.\ Cor.\ \ref{CorKM})
the number of deformation types of those is finite.\qed

\bigskip

Note that the argument does not show the finiteness of the number of connected
components of $\gM_\Gamma$ that parametrizes complex structures on a fixed
manifold $M$ is finite. A priori, this seems only to be the case modulo
the action of $\Aut(\Gamma)$.

The theorem is equivalent to a statement about the number of connected
components of the space of all $\SSp(n)$-metrics on $M$. Since the pull-back
of an $\SSp(n)$-metric under a diffeomorphism of $M$ yields again an
$\SSp(n)$-metric on $M$, one only gets a finiteness result modulo the action
of the diffeomorphism group.

\begin{theorem}---
The group ${\rm Diff}(M)/{\rm Diff}^0(M)$ acts on the set
of connected components of the set of all $\SSp(n)$-metrics on $M$.
The quotient by this action is finite.\qed
\end{theorem}

\section{Using the Hitchin-Sawon formula}

Applying the theory of Rozansky-Witten invariants, Hitchin and Sawon in
\cite{HS} proved the following formula

\begin{theorem}{\bf \cite[Thm.\ 3]{HS}} ---
Let $(M,g)$ be a compact Riemannian manifold of dimension $4n$ with holonomy
$\SSp(n)$. Then
$$\frac{||R||^{2n}}{(192\pi^2n)^n}=\int_M\sqrt{\hat A(M)}\cdot(\vol M)^{n-1}.$$
\end{theorem}

Here, $R$ is the curvature of $(M,g)$ and $\hat A$ is the
$\hat A$-genus of $M$. By fixing a complex structure $I$ on $M$ that is
compatible
with the hyperk\"ahler metric $g$ we obtain an irreducible holomorphic
symplectic manifold $X=(M,I)$. Then $g$ and $I$ define a K\"ahler form
$\omega_I$ on $X$. With respect to this K\"ahler form the norm of the curvature
can be expressed as (cf.\ \cite{Be})
$||R||^2=\frac{8\pi^2}{(2n-2)!}\int_Xc_2(X)\omega_I^{2n-2}$ and for the
volume of $M$ one finds $\vol(M)=(1/(2n)!)\int_X\omega_I^{2n}$.
Furthermore, $\hat A(M)=\td(X)$.
Hence, the Hitchin-Sawon formula can be rewritten as
$$\left(\int_X c_2\omega_I^{2n-2}\right)^n=c\int_X\sqrt{\td(X)}\cdot
\left(\int_X\omega_I^{2n}\right)^{n-1}\eqno{(2)}$$
with $c=\frac{(24n)^n}{((2n)!)^{n-1}}$. As is pointed out in
\cite{HS} this immediately yields
$$\int_X\sqrt{\td(X)}=\int_M\sqrt{\hat A(M)}>0.$$
Clearly, $c\int_X\sqrt{\td(X)}$ is a rational number which can be
written as $p/q$ with positive integers $p$ and $q$, where $q$ is
bounded from above by a universal constant $c_n$ depending only on $n$.
Hence, $c\int_X\sqrt{\td(X)}\geq1/c_n$. Since
any K\"ahler class $\omega$ on an irreducible holomorphic symplectic manifold
$X$ is obtained as $[\omega_I]$, this yields

\begin{corollary}---
Let $X$ be an irreducible holomorphic symplectic manifold of dimension
$2n$. If $\omega$ is any K\"ahler class on $X$, then $\int_X\omega^{2n}<
c_n^{1/(n-1)}\cdot(\int_Xc_2(X)\omega^{2n-2})^{\frac{n}{n-1}}$, where $c_n$ depends only
on $n$.\qed
\end{corollary}

\begin{theorem}\label{mainthm}---
\it If the second integral cohomology $H^2$
and the homogeneous polynomial of degree $2n-2$ on $H^2$ defined
by the first Pontrjagin class are given, then there exist at most
finitely many diffeomorphism types of compact
hyperk\"ahler manifolds of dimension $4n$ realizing this structure.
\end{theorem}

\prf Let $M$ be any compact real $4n$-dimensional manifold and
let $I$ be a complex structure on $M$. Then, $p_1(M)=-2c_2(X)$, where
$X=(M,I)$. In particular, the action of $c_2(X)$ on
$H^2(X,\IZ)=H^2(M,\IZ)$ does not depend on the chosen complex structure $I$.
By assumption the second cohomology $H^2$ of $M$ as an abelian
group and the homogeneous polynomial  $\alpha\mapsto\int_X\alpha^{2n-2}p_1$
of degree $2n-2$ on $H^2$ are fixed. If $X=(M,I)$ is an irreducible
holomorphic symplectic manifold that admits an ample line bundle $L\in\Pic(X)$,
then $L^{2n}$ is bounded from above by
$c_n^{1/(n-1)}(\int_X \frac{-p_1}{2}.L^{2n-2})^\frac{n}{n-1}$. Thus,
by Koll\'ar-Matsusaka the number of deformation types of these complex
manifolds is finite. Using the surjectivity of the period map and arguing
as in the proof of Thm.\ \ref{mainthm2}, one concludes that the number of
diffeomorphism types of compact manifolds $M$ of dimension $4n$
carrying a Riemannian metric with holonomy $\SSp(n)$ is finite.\qed

\bigskip

\begin{corollary}---
The number of deformation types of irreducible holomorphic symplectic manifolds
$X$ of dimension $2n$ with given $H^2(X,\IZ)$ (as abelian group) and
given homogeneous polynomial $H^2(X,\IZ)\to\IZ$,
$\alpha\mapsto\int_X\alpha^{2n-2}c_2(X)$
is finite.\qed
\end{corollary}

\section{Finiteness for fixed quadratic form}\label{unnorm}

Let $X$ be an irreducible holomorphic symplectic manifold and
$q_X$ the Beauville-Bogomolov form on $H^2(X,\IZ)$.
If $\sigma$ generates $H^0(X,\Omega^2_X)$ such that
$\int_X(\sigma\bar\sigma)^n=1$, then $q_X$ is a scalar multiple
of the form
$f_X(\alpha)=(n/2)\int_X\beta^2(\sigma\bar\sigma)^{n-1}+\lambda\mu$, where
$\alpha=\lambda\sigma+\mu\bar\sigma+\beta$ is the Hodge-decomposition of
$\alpha$.
Usually, $q_X$ is chosen such that it becomes a primitive integral form on
$H^2(X,\IZ)$. Note that $q_X$ usually does not define the structure of a
unimodular lattice on $H^2(X,\IZ)$. So classification theory of unimodular
lattices cannot be applied.
Due to a result of Fujiki (cf.\ \cite[1.11]{H} or \cite{Bo3}) one has

\begin{theorem}{\bf \cite[4.12]{Fujiki2}} ---
There exists a constant $c\in\IQ$ such that for all $\alpha\in H^2(X)$
$$c\cdot q_X(\alpha)=\int_X\sqrt{\td(X)}\alpha^2.$$
\end{theorem}

Note that $\sqrt{\td(X)}$ does not depend on the complex structure $I$. Since
all odd Chern classes of $X$ are trivial, $\sqrt{\td(X)}$ is a universal
expression in the Pontrjagin classes of the underlying real manifold $M$.

A priori, the constant $c$ could be trivial or negative. That this is not the
case follows from a result of Nieper \cite{N}

\begin{theorem}---
For any $\alpha\in H^2(X)$ one has $$\int_X\sqrt{\td(X)}{\rm
exp}(\alpha)=(1+\lambda(\alpha))^n\int_X\sqrt{\td(X)},$$ where $\lambda$ is a
positive multiple of the quadratic form $f_X$ (or equivalently of $q_X$).
\end{theorem}

In fact $\lambda$ and $f_X$ differ by
$12/\int_Xc_2(X)(\sigma\bar\sigma)^{n-1}$. The formula in particular gives back
(1). More generally, it yields $$\int_X\sqrt{\td(X)}\alpha^{2k}=(2k)!{n\choose
k}\lambda(\alpha)^k\int_X\sqrt{\td(X)}.$$ For $k=1$ one obtains
$\int_X\sqrt{\td(X)}\alpha^2=2{n\choose 2}\int_X\sqrt{\td(X)}\lambda(\alpha)$.
Thus the constant $c$ in Fujiki's result is a positive multiple of
$(\int_X\sqrt{\td(X)})^{-1}$. The latter is positive as was observed by Hitchin
and Sawon. Note that Nieper's result also nicely generalizes the Hitchin-Sawon
formula (2) to
$$\left(\int_X\sqrt{\td(X)}\alpha^{2k}\right)^n=\left((2k)!{n\choose
k}\right)^n
\left(\frac{1}{(2n)!}\int_X\alpha^{2n}\right)^k\cdot\left(\int_X\sqrt{\td(X)}\right)^{n-k}.$$

Hence, for two different irreducible holomorphic symplectic structures
$X=(M,I)$ and $X'=(M,I')$ on a manifold $M$ the associated quadratic
forms $q_X$ and $q_{X'}$ differ by a positive multiple and, therefore,
are equal. This was used in Sect.\ \ref{complex}.

Let us define the {\it unnormalized Beauville-Bogomolov form} $\tilde q_X$ on
$H^2(X)$ by $\tilde q_X(\alpha):=d_n\int_X\sqrt{\td(X)}\alpha^2$, where $d_n$ is
a universal integer (of no importance) depending only on $n$ such that
$d_n(\sqrt{\td(X)})_{4n-4}$ is universally an integral class. This quadratic
form $\tilde q_X$ is a positive integral multiple of the primitive form $q_X$.
As one usually does not have much information about the divisibility of Chern
classes, it looks, at least from this point of view, more naturally to work
with $\tilde q_X$ instead of the primitive form $q_X$.

As it is widely conjectured, a global Torelli theorem for higher dimensional
hyperk\"ahler manifold should assert that two irreducible holomorphic
symplectic manifolds are birational whenever there exists an isomorphism
$H^2(X,\IZ)\cong H^2(X',\IZ)$ that respects Hodge structures and
the Beauville-Bogomolov quadratic forms $q_X$ and $q_{X'}$, resp.
A modified version of it would replace $q_X$ and $q_{X'}$ by the unnormalized
Beauville-Bogomolov forms $\tilde q_X$ and $\tilde q_{X'}$, resp. Note
that the modified version is a consequence of the original one, but it
is, a priori, weaker. Here one uses that a birational map
$X  - \to X'$ induces an
isomorphism $H^*(X)\cong H^*(X')$ that maps $c(X)$ to $c(X')$.

From this point of view, it seems reasonable to
replace $(H^2(X,\IZ),q_X)$ by $(H^2(X,\IZ),\tilde q_X)$ also in the finiteness
question we are interested here. The proof of the following theorem
is omitted as it uses arguments already explained in the
previous sections.

\begin{theorem} ---
There exists only a finite number of deformation types of irreducible
holomorphic symplectic manifolds $X$ of fixed dimension such that
the lattice $(H^2(X,\IZ),\tilde q_X)$ is isomorphic to a given one.
\end{theorem}

It should be clear from the discussion that fixing the lattice
$(H^2(X,\IZ),\tilde q_X)$ is actually too much.
Everything that is needed is a class $\alpha\in H^2(X,\IZ)$ with
fixed $\tilde q_X(\alpha)>0$. Thus, for any $k\in\IZ_{>0}$ there is
only a finite number of deformation types of irreducible
holomorphic symplectic manifolds $X$ that admit an element $\alpha\in
H^2(X,\IZ)$ with $\tilde q_X(\alpha)=k$.

\bigskip

{\bf Acknowledgement:} I would like to thank
Manfred Lehn, Marc Nieper, and Justin Sawon
for comments and useful conversations.
\bigskip

{\footnotesize 

\bigskip

\bigskip

\noindent
{\small Mathematisches Institut\\
Universit\"at zu K\"oln\\
Weyertal 86-90\\
50931 K\"oln, Germany\\
\texttt{huybrech@mi.uni-koeln.de}}

\end{document}